\documentclass[12pt,a4paper]{amsart}
\usepackage[utf8]{inputenc}	
\usepackage[T1]{fontenc}
\usepackage[in,headings]{fullpage}
\usepackage{amsbsy, amscd, amsfonts, amsmath, amsrefs, amssymb, amsthm}
\usepackage{graphicx}
\usepackage{float}
\usepackage{color}   
\usepackage[colorlinks=true,linkcolor=blue,citecolor=blue,urlcolor=blue]{hyperref}

\newtheorem{theorem}{Theorem}

\newtheorem{corollary}[theorem]{Corollary}

\newtheorem{lemma}[theorem]{Lemma}

\theoremstyle{definition}

\newtheorem{remark}[theorem]{Remark}

\theoremstyle{remark}

\newcommand{\C}{\mathbf{C}}

\newcommand{\R}{\mathbf{R}}

\renewcommand{\Re}{\mathop{\mathrm{Re}}\nolimits}
\renewcommand{\Im}{\mathop{\mathrm{Im}}\nolimits}
\newcommand{\Rzeta}{\mathop{\mathcal R }\nolimits}

\newfont{\cmbsy}{cmbsy10}
\newfont{\cmmib}{cmmib10}
\newcommand{\Orden}{\mathop{\hbox{\cmbsy O}}\nolimits}
\newcommand{\orden}{\mathop{\hbox{\cmmib o}}\nolimits}

\begin{document}

\title[Region without zeros for $\Rzeta(s)$]{Regions without zeros for the auxiliary function of Riemann. }
\author[Arias de Reyna]{J. Arias de Reyna}
\address{%
Universidad de Sevilla \\ 
Facultad de Matem\'aticas \\ 
c/Tarfia, sn \\ 
41012-Sevilla \\ 
Spain.} 




\email{arias@us.es, ariasdereyna1947@gmail.com}


\begin{abstract}
We give explicit and extended versions of some of Siegel's results. We extend the validity of Siegel's asymptotic development in the second quadrant to most of the third quadrant. We also give precise bounds of the error; this allows us to give an explicit region free of zeros, or with only trivial zeros.  The left limit of the zeros on the upper half plane is extended from  $1-\sigma\ge a t^{3/7}$ in Siegel to $1-\sigma\ge A t^{2/5}\log t$.  Siegel claims that it can be proved that there are no zeros in the region $1-\sigma\ge t^\varepsilon$ for any $\varepsilon>0$. We show that Siegel's proof for the exponent $3/7$ does not extend to prove his claim.
\end{abstract}

\maketitle

\section{Introduction}

Siegel introduced in \cite{Siegel} the function 
\begin{equation}
\Rzeta(s)=\int_{0\swarrow1}\frac{x^{-s} e^{\pi i x^2}}{e^{\pi i x}-
e^{-\pi i x}}\,dx,
\end{equation}
that Riemann connected with his function $\zeta(s)$. Siegel says that Riemann tried to obtain information concerning the zeros of $\Rzeta(s)$ but that very little is known about them. 

Here, refining Siegel's results,
we give explicit regions free of zeros of Riemann's auxiliar function.
In order to do this, we give concrete form to Siegel's theorems in \cite{Siegel}. In particular, we study the left limit of the zeros $\rho=\beta+i\gamma$ of $\Rzeta(s)$ with $\gamma>0$.  Also, we prove the result of Siegel on the asymptotic behavior 
of $\Rzeta(s)$ on the second quadrant extending them to most of the third quadrant.

There are some differences with Siegel. First, we extend the region where the asymptotic representation is valid. Extending the validity to almost all the region 
$\sigma<0$, (see precise formulation in Theorems \ref{mainTheorem} and \ref{apRzeta}). We also give precise bounds of the error; this allows us to give an explicit region free of zeros, or with only trivial zeros. 

There is also a note about Siegel's paper. He claims that he can prove that for every $\varepsilon>0$ there exists $t_0$ such that the region $t>t_0$ and $1-\sigma\ge c t^{\varepsilon}$ does not contain zeros of $\Rzeta(s)$. Siegel only proves the case $1-\sigma\ge t^{3/7}$.  I try to follow Siegel's argument in this case, but I have been unable to obtain Siegel's claim. My arguments give only the extension to  $1-\sigma\ge t^{2/5}\log t$. It remains  an open problem to find the true limit of the zeros. 
My computation of zeros \cite{A172} shows that there are zeros to the left of the critical line and very far from it.

Another revision of Siegel's paper has been posted recently \cite{O'S}, which do not address this problem.

\section{Application of the Cauchy Theorem.}

\begin{theorem}\label{mainTh}
Let $s=\sigma+it$ with $0<\arg(s-1)<2\pi$, we define
\begin{gather}
\eta:=\sqrt{\frac{s-1}{2\pi i}},\quad-\frac{\pi}{4}<\arg\eta<\frac{3\pi}{4},\quad\eta_1:=\Re(\eta),\quad
\eta_2:=\Im(\eta),\label{E:etadef}\\ m:=\lfloor\eta_1+\eta_2\rfloor,\quad
w(z):=\exp\Bigl\{2\pi i\eta^2\Bigl(\log\Bigl(1+\frac{z}{\eta}\Bigr)-
\frac{z}{\eta}+\frac12\Bigl(\frac{z}{\eta}\Bigr)^2\Bigr)\Bigr\}-1.\label{defw}
\end{gather}
For each integer $k\ge m$ we have
\begin{multline}\label{main}
\Rzeta(s)=\chi(s)\Bigl\{\zeta(1-s)-\sum_{n=1}^{k}n^{s-1}-\\-\eta^{s-1}e^{-\pi i \eta^2}
\Bigl[\frac{\sqrt{2}e^{3\pi i/8}\sin\pi\eta-(-1)^k e^{2\pi i\eta-2\pi i(\eta-k)^2}}
{2\cos2\pi\eta}+\\+\sum_{j=m+1}^k(-1)^{j-1}e^{-2\pi i(j-\eta)^2}
w(j-\eta)+R\Bigr]\Bigr\}
\end{multline}
where  $R$ is given by
\begin{equation}\label{defR}
R:=\int_{m\searrow m+1}\frac{e^{-2\pi i(x-\eta)^2}}{e^{\pi i x}-e^{-\pi i
x}}w(x-\eta)\,dx.
\end{equation}
\end{theorem}

\begin{remark}
Compare Siegel \cite{Siegel}*{(65), (66), (67),  and (68)}.
\end{remark}

\begin{proof}
We start with the equation (for definition and basic results about the Riemann auxiliary function $\Rzeta(s)$ see \cite{A166})
\begin{equation}
\Rzeta(s)=\chi(s)\bigl\{\zeta(1-s)-\overline{\Rzeta(1-\overline{s})}\bigr\}.
\end{equation}
and observe that 
\begin{displaymath}
\overline{\Rzeta(1-\overline{s})}=\overline{
\int_{0\swarrow1}\frac{x^{\overline{s}-1}e^{\pi i x^2}}{e^{\pi i x}-e^{-\pi i
x}}
\,dx}=\int_{0\searrow1}
\frac{x^{s-1}e^{-\pi i x^2}}{e^{\pi i x}-e^{-\pi i x}}\,dx.
\end{displaymath}
Therefore, we have 
\begin{equation}
\Rzeta(s)=\chi(s)\Bigl\{\zeta(1-s)-
\int_{0\searrow1}
\frac{x^{s-1}e^{-\pi i x^2}}{e^{\pi i x}-e^{-\pi i x}}\,dx\Bigr\}.
\end{equation}

The saddle-point $\eta$ of the function $x^{s-1}e^{-\pi i x^2}$ is $\eta=\sqrt{\frac{s-1}
{2\pi i}}$. Since we asume $0<\arg(s-1)<2\pi$ we may take $-\pi/4<\arg\eta<\frac{3\pi}{4}$, fixing in this way the value of $\eta$. We will denote by $\eta_1$ and $\eta_2$ the real and imaginary part of $\eta=\eta_1+i\eta_2$. Let $k\ge0$ be an integer, by Cauchy's Theorem we have
\begin{displaymath}
\int_{0\searrow1}
\frac{x^{s-1}e^{-\pi i x^2}}{e^{\pi i x}-e^{-\pi i x}}\,dx=
\sum_{n=1}^kn^{s-1}+\int_{k\searrow k+1}
\frac{x^{s-1}e^{-\pi i x^2}}{e^{\pi i x}-e^{-\pi i x}}\,dx.
\end{displaymath}
Of course, this will be more convenient when the new path of integration contains the saddle-point $\eta$ or is near it, that is, when $k=m$. But following Siegel we use here a general $k$. The restriction on $s$ and the election of $\eta$  is given so that $m\ge0$. Therefore,  we may move the line of integration to $m\searrow m+1$ on the region where the integrand is meromorphic ($-\pi<\arg x <\pi$).

 The integrand may be transformed as usual (see \cite{A86}, or \cite{G})
\begin{multline}
x^{s-1}e^{-\pi i x^2}=e^{2\pi i \eta^2 \log(\eta+x-\eta)}e^{-\pi i (x-\eta+\eta)^2}
\\=\eta^{s-1}e^{-\pi i \eta^2}\exp\Bigl\{2\pi i \eta^2
\Bigl(\log\Bigl(1+\frac{x-\eta}{\eta}\Bigr)-\frac{x-\eta}{\eta}+\frac{(x-\eta)^2}{2\eta^2}\Bigr)\Bigr\}e^{-2\pi i (x-\eta)^2}\\
=\eta^{s-1}e^{-\pi i \eta^2}(w(x-\eta)+1)e^{-2\pi i (x-\eta)^2}.
\end{multline}
where  $w(z)$ is defined in \eqref{defw}.

It follows that 
\begin{multline}
\int_{0\searrow1}
\frac{x^{s-1}e^{-\pi i x^2}}{e^{\pi i x}-e^{-\pi i x}}\,dx
=\sum_{n=1}^kn^{s-1}+\\
+\eta^{s-1}e^{-\pi i \eta^2}\Bigl\{
\int_{k\searrow k+1}\frac{e^{-2\pi i(x-\eta)^2}}{e^{\pi i x}-e^{-\pi i x}}\,dx
+
\int_{k\searrow k+1}\frac{e^{-2\pi i(x-\eta)^2}}{e^{\pi i x}-e^{-\pi i
x}}w(x-\eta)\,dx\Bigr\}.
\end{multline}

The first integral may be computed explicitly (see \cite{Siegel}*{eq.~(67)}, this is due in fact to  Riemann) 
\begin{equation}
\int_{k\searrow k+1}\frac{e^{-2\pi i(x-\eta)^2}}{e^{\pi i x}-e^{-\pi i x}}\,dx
=\frac{\sqrt{2}e^{3\pi i/8}\sin\pi\eta-(-1)^k e^{2\pi i\eta-2\pi i(\eta-k)^2}}
{2\cos2\pi\eta}.
\end{equation}

Now we may apply again Cauchy's Theorem, for example since we assume  that $k\ge m$ we get
\begin{multline}
\int_{k\searrow k+1}\frac{e^{-2\pi i(x-\eta)^2}}{e^{\pi i x}-e^{-\pi i
x}}w(x-\eta)\,dx\\
=\sum_{j=m+1}^k(-1)^{j-1}e^{-2\pi i(j-\eta)^2}w(j-\eta)+
\int_{m\searrow m+1}\frac{e^{-2\pi i(x-\eta)^2}}{e^{\pi i x}-e^{-\pi i
x}}w(x-\eta)\,dx.
\end{multline}
\end{proof}

\section{Bound of the remainder.}

What follows is directed to bounding $R$ in the above Theorem. 

\begin{lemma}\label{lemainteg}
Let $w=u+iv\in\C$ with $u$ and $v\in\R$. If $u+|v|\ge-1/2$ then 
\begin{equation}
\Bigl|\int_0^w\frac{t^2}{1+t}\,dt\Bigr|\le \frac{3}{4}|w|^2
\end{equation}
integrating along the segment joining $0$ and $w$. 
\end{lemma}

\begin{proof}
The integral defines a holomorphic function  $f(w)$ on the plane with a cut along the negative real axis from $-1$ to $-\infty$.  Condition $u+|v|\ge-1/2$ excludes an angle that contains this cut.  

The proof is somewhat complicated, since it is not obtained integrating the absolute value of the integrand.

Since $f(w)$ is real for real $w$, we have $f(\overline{w})=\overline{f(w)}$ and we can assume that $v\ge0$. 

Put $b=(1+i)/2$ and $T(w)=w/(w+b)$. We have $T(0)=0$, $T(-b)=\infty$, since $-b$ is the conjugate of $0$ with respect to the line $x+y=-1/2$, $T$ transform this line on a circle with center at $0$. $T(-1/2)=i$, so that this circle is the unit circle.
Our inequality may be written as 
\begin{displaymath}
\Bigl|\frac{1}{(w+b)^2}\int_0^w\frac{t^2}{1+t}\,dt\Bigr|\le\frac{3}{4}\Bigl|\frac{w}
{w+b}\Bigr|^2.
\end{displaymath}
For the values of $w$ we are interested in, we have $|w/(w+b)|=|T(w)|\le 1$. 

Our integral may be computed explicitly.
\begin{equation}
\int_0^w\frac{t^2}{1+t}\,dt=\log(1+w)-w+\frac{1}{2}w^2.
\end{equation}
Changing variables $z$ instead of $w$ with
\begin{displaymath}
z=\frac{w}{w+b},\quad w=\frac{bz}{1-z}, \quad w+b=\frac{b}{1-z}
\end{displaymath}
our inequality may be written as 
\begin{equation}
\Bigl|\Bigl(\frac{1-z}{b}\Bigr)^2\log\frac{1-z(1-b)}{1-z}-\frac{z}{ b}+\frac{z^2}{
b}+\frac{z^2}{2}\Bigr|\le\frac{3}{4}|z|^2,\qquad |z|<1.
\end{equation}
We expand  in Taylor series the function to be bounded.
\begin{displaymath}
\Bigl(\frac{1-z}{ b}\Bigr)^2\Bigl(\sum_{n=1}^\infty\frac{z^n}{ n}-\sum_{n=1}^\infty
\frac{(1-b)^nz^n}{ n}\Bigr)-\frac{z}{b}+\frac{z^2}{
b}+\frac{z^2}{2}.
\end{displaymath}
The first term cancels out, so that we must prove the inequality
\begin{equation}
\Bigl|\sum_{n=3}^\infty\frac{z^n}{b^2}\Bigl(\frac{1}{ n}-\frac{2}{ n-1}+\frac{1}{n-2}
-\frac{(1-b)^n}{n}+\frac{2(1-b)^{n-1}}{n-1}-\frac{(1-b)^{n-2}}{ n-2}\Bigr)
\Bigr|\le\frac{3}{4}|z|^2.
\end{equation}
The power series here has convergence radius $=1$ and converges absolutely for all $|z|\le1$. 
By Schwarz's lemma the proof reduces to check the inequality
\begin{equation}
\sum_{n=3}^\infty\Bigl|\frac{1}{b^2}\Bigl(\frac{1}{n}-\frac{2}{n-1}+\frac{1}{n-2}
-\frac{(1-b)^n}{n}+\frac{2(1-b)^{n-1}}{n-1}-\frac{(1-b)^{n-2}}{n-2}\Bigr)
\Bigr|\le\frac{3}{4}.
\end{equation}
The sum is equal to $0.7439893\dots$ so that this inequality is true.

\end{proof}

\begin{lemma}\label{lemmaboundcos}
For all real $u$, we have $|\cos((1-i)u)|^2>\frac14 \cosh(2u)$, so that 
\begin{equation}\label{E:cose}
|\cos((1-i)u)|> 2^{-3/2} e^{|u|}.
\end{equation}
\end{lemma}

\begin{proof}
For real $u$, we have $|\cos((1-i)u)|^2=\frac12(\cos(2u)+\cosh(2u))$. The first inequality is then equivalent to $\cosh(2u)\ge-2\cos(2u)$.  Then \eqref{E:cose} follows easily. 

\end{proof}

\begin{theorem}\label{boundR}
For $\sigma<1$ and $|s-1|\ge4 \pi$ we have \eqref{main} with $R$ as in \eqref{defR} bounded by
\begin{equation}\label{boundReq1}
|R|\le 
\frac{\min(4, 29e^{-\pi \eta_2})}{|\eta|}+ 15\frac{e^{-\frac{\pi}{32}|\eta|^{2}}}{|\eta|}.
\end{equation}
\end{theorem}

\begin{proof}
$R$ is given by the integral in \eqref{defR}. The path of integration  is a straight line $L$ of  direction $e^{-\pi i/4}$ through the point $m+1/2$. The point $\eta$ is contained in the strip of points with a distance $\le 2^{-3/2}$ to this line. Let $\eta'$ be the point of the line of integration nearest to $\eta$, then $\eta'-\eta=\alpha e^{\pi i/4}$ with $\alpha\in\R$ and  $|\alpha|\le 2^{-3/2}$ and the line $L$ is parametrized by $x=\eta'+ue^{-\pi i/4}$ with $u\in\R$. 

We are going to bound each term of the integrand in \eqref{defR}. 
First,
\begin{multline*}
-2\pi i(x-\eta)^2=-2\pi i(\eta'-\eta+ue^{-\pi i/4})^2=-2\pi i(\alpha e^{\pi i/4}+
ue^{-\pi i/4})^2=\\=2\pi\alpha^2-2\pi u^2-4\pi \alpha u i
\end{multline*}
so that 
\begin{equation}
|e^{-2\pi i(x-\eta)^2}|=e^{2\pi (\alpha^2-u^2)}\le e^{\pi/4}e^{-2\pi u^2}.
\end{equation}

Since $\eta'$ is on the line $L$ we can put it in the form $\eta'=m+1/2-re^{-\pi i /4}$ for some real $r$, if we put $\eta'=\eta'_1+i\eta'_2$, we will have 
$\eta'_2=r/\sqrt{2}$, and from 
$\eta'=\alpha e^{\pi i/4}+\eta$ we will get $\eta'_2=\eta_2+\alpha/\sqrt{2}$, so that
\begin{equation}\label{valuer}
r=\alpha+\sqrt{2}\,\eta_2.
\end{equation}
Then $x=\eta'+ue^{-\pi i/4}=m+1/2+(u-r)e^{-\pi i/4}$ and
\begin{multline*}
e^{\pi i x}-e^{-\pi i x}=2i\sin\pi x=2i\sin\bigl(\pi (m+1/2+(u-r)e^{-\pi i/4})\bigr)=\\=2i(-1)^m \cos\Bigl(\frac{\pi(u-r)}{\sqrt{2}}(1-i)\Bigr)
\end{multline*}
by Lemma \ref{lemmaboundcos} we get
\begin{equation} \label{boundcos}
|e^{\pi i x}-e^{-\pi i x}|\ge 2^{-1/2}e^{\pi\frac{|r-u|}{\sqrt{2}}}\ge 
2^{-1/2}\max(1, e^{\frac{\pi(r-u)}{\sqrt{2}}}).
\end{equation}

We want to apply Lemma \ref{lemainteg} to $w=\eta^{-1}(x-\eta)$ for $x$ on the integration path. We have seen that $x-\eta=\alpha e^{\pi i/4}+ue^{-\pi i/4}$ with $\alpha$, $u\in\R$ and $|\alpha|\le 2^{-3/2}$. Hence $x-\eta$ is contained in a strip of direction $e^{-\pi i /4}$ and width $2|\alpha|$ centered at the origin. Set $\eta=|\eta|e^{i\varphi}$. The hypothesis $\sigma\le 1$ implies $0<\varphi<\frac{\pi}{2}$. It follows that $\frac{x-\eta}{\eta}=|\eta|^{-1}(\alpha e^{i(\frac{\pi}{4}-\varphi)}+ue^{-i(\frac{\pi}{4}+\varphi)})$ is contained in a strip of direction $e^{-i(\frac{\pi}{4}+\varphi)}$ and width $2|\alpha|/|\eta|$ centered at the origin. We need this strip to be contained in the region  $\Omega:=\{u+iv\in\C: u+|v|>-1/2\}$ (~the region where Lemma \ref{lemainteg} applies~). For this, we need the direction to be in the range  $-\frac{3\pi}{4}<\frac{\pi}{4}+\varphi<\frac{3\pi}{4}$. This is true by our hypothesis $\sigma<1$. The strip will be in $\Omega$ if the intersection of the strip with the real axis is contained in $(-1,+\infty)$. But the leftmost point of this intersection is 
$-|\alpha|/(|\eta|\sin(\frac{\pi}{4}+\varphi))$. The needed inequality is  true,
since for $0<\varphi<\frac{\pi}{2}$ and  $|\eta|>1$ we have
\[\frac{|\alpha|}{\sin(\varphi+\frac{\pi}{4})}\le \frac{2^{-3/2}}{2^{-1/2}}=\frac{1}{2}\le \frac{|\eta|}{2}.\]

Therefore, by Lemma \ref{lemainteg} we have
\begin{displaymath}
\Bigl|\log\Bigl(1+\frac{x-\eta}{\eta}\Bigr)-\frac{x-\eta}{\eta}+\frac12\Bigl(\frac{x-\eta}{\eta}\Bigr)^2\Bigr|=\Bigl|\int_0^{
\frac{x-\eta}{\eta}}\frac{t^2\,dt}{1+t}\Bigr|\le \frac{3}{4}\frac{|x-\eta|^2}{|\eta|^2}.
\end{displaymath}
Since $|e^z-1|\le e^{|z|}-1\le e^{|z|}$
and by  definition \eqref{defw} of  $w(z)$  for $x\in L$ 
\begin{equation}\label{generalineq}
|w(x-\eta)|\le e^{\frac{3\pi}{2}|x-\eta|^2}=e^{\frac{3\pi}{2}(\alpha^2+u^2)}
\le e^{\frac{3\pi}{16}}e^{\frac{3\pi}{2}u^2}.
\end{equation}

For the points on the strip $S$ we may obtain another bound applicable when 
$|x-\eta|<\frac12|\eta|$. In fact, for $|z|<\frac12$ we have
\begin{displaymath}
|\log(1+z)-z+\tfrac12 z^2|=\Bigl|\sum_{n=3}^\infty \frac{(-1)^{n-1}}{n}z^n\Bigr|\le 
\frac13\frac{|z|^3}{1-|z|}\le \frac23 |z|^3
\end{displaymath}
so that for $|x-\eta|\le\frac12|\eta|$ 
\begin{equation}\label{restrineq}
|w(x-\eta)|\le e^{\frac{4\pi}{3}\frac{|x-\eta|^3}{|\eta|}}-1.
\end{equation}
We divide the integral \eqref{defR} for $R$ into two, one for $|u|<\frac14|\eta|$ and the rest. In the first case $|x-\eta'|=|u|$ and  
\begin{displaymath}
|x-\eta|\le |x-\eta'|+|\eta'-\eta|\le \frac14|\eta|+2^{-3/2}\le
\tfrac12|\eta|
\end{displaymath}
(we apply here that from  $|s-1|\ge4\pi$ we get  $2^{-3/2}\le\frac14|\eta|$).  So in the range
$|u|<\frac14|\eta|$ we may apply 
\eqref{restrineq} and for the rest  $u$ we may apply 
\eqref{generalineq}. Therefore,
\begin{multline}\label{split}
|R|\le \int_{|u|<|\eta|/4}\frac{e^{\pi/4} e^{-2\pi u^2}}{2^{-1/2} 
e^{\pi(r-u)/\sqrt{2}}}\;\bigl(e^{\frac{4\pi}{3}\frac{|x-\eta|^3}{|\eta|}}-1\bigr)\,du+\\ +
\int_{|u|>|\eta|/4}\frac{e^{\pi/4} e^{-2\pi u^2}}{2^{-1/2} 
}e^{3\pi/16}e^{\frac{3\pi}{2}u^2}\,du\le 
\end{multline}
\begin{multline*}
\le \sqrt{2} e^{\pi/4}e^{-\pi r/\sqrt{2}}\int_{|u|<|\eta|/4}
e^{-2\pi u^2}\bigl(e^{\frac{4\pi}{3}\frac{(\alpha^2+u^2)^{3/2}}{|\eta|}}-1\bigr)
e^{\pi u/\sqrt{2}}\,du+\\+2\sqrt{2} e^{7\pi/16}\frac{4}{|\eta|}
\int_{|\eta|/4}^{+\infty}u e^{-\pi u^2/2}\,du:= I_1+I_2.
\end{multline*}
We may compute exactly $I_2$ 
\begin{displaymath}
I_2=\frac{8\sqrt{2}e^{7\pi/16}}{\pi}\frac{e^{-\pi|\eta|^2/32}}{|\eta|}<14.2355
\frac{e^{-\pi|\eta|^2/32}}{|\eta|}.
\end{displaymath}
By \eqref{valuer} we have $e^{-\pi r/\sqrt{2}}=e^{-\pi\alpha/\sqrt{2}}e^{-\pi\eta_2}
\le e^{\pi/4}e^{-\pi\eta_2}$. Then 
\begin{multline*}
I_1<2\sqrt{2} e^{\pi/2}e^{-\pi \eta_2}\int_0^{|\eta|/4}
e^{-2\pi u^2}\bigl(e^{\frac{4\pi}{3}\frac{(\alpha^2+u^2)^{3/2}}{|\eta|}}-1\bigr)
e^{\pi u/\sqrt{2}}\,du
<\\
<2\sqrt{2} e^{\pi/2}e^{-\pi \eta_2}\int_0^{|\alpha|}
e^{-2\pi u^2+\pi u/\sqrt{2}}\bigl(e^{\frac{4\pi}{3}\frac{2^{3/2}|\alpha|^3}{|\eta|}}-1\bigr)\,du+\\
+2\sqrt{2} e^{\pi/2}e^{-\pi \eta_2}\int_0^{|\eta|/4}
e^{-2\pi u^2+\pi u/\sqrt{2}}\bigl(e^{\frac{4\pi}{3}\frac{2^{3/2}u^3}{|\eta|}}-1\bigr)\,du:=I_{11}+I_{12}.
\end{multline*}
For $0\le x\le 1$ we have $|e^x-1|\le2x$, so that  since $|\alpha|<2^{-3/2}$ 
and $|\eta|\ge\sqrt{2}$ 
\begin{displaymath}
\bigl|e^{\frac{4\pi}{3}\frac{2^{3/2}|\alpha|^3}{|\eta|}}-1\bigr|\le 2
\frac{4\pi}{3}\frac{2^{3/2}|\alpha|^3}{|\eta|}\le \frac{\pi}{3|\eta|}.
\end{displaymath}
It follows that 
\begin{displaymath}
I_{11}\le \frac{2\pi\sqrt{2} e^{\pi/2}}{3}\frac{e^{-\pi \eta_2}}{|\eta|}
\int_{0}^{2^{-3/2}}e^{-2\pi u^2+\pi u/\sqrt{2}}\,du<5.7518 \frac{e^{-\pi \eta_2}}{|\eta|}.
\end{displaymath}

Furthermore, since $\frac{e^x-1}{x}<e^x$ for $x>0$
\begin{displaymath}
I_{12}<2\sqrt{2} e^{\pi/2}e^{-\pi \eta_2} \frac{4\pi}{3}\frac{2^{3/2}}{|\eta|}\int_0^{|\eta|/4}
e^{-2\pi u^2+\pi u/\sqrt{2}}e^{\frac{4\pi}{3}\frac{2^{3/2}u^3}{|\eta|}} u^3\,du.
\end{displaymath}
By the dominated convergence Theorem, the last integral, for $|\eta|\to+\infty$,  has a  finite limit 
\begin{displaymath}
\int_0^{+\infty}u^3 e^{-2\pi u^2+\pi u/\sqrt{2}}\,du=0.0453198
\end{displaymath}
a numerical study shows that 
\begin{displaymath}
\int_0^{|\eta|/4}
e^{-2\pi u^2+\pi u/\sqrt{2}}e^{\frac{4\pi}{3}\frac{2^{3/2}u^3}{|\eta|}} u^3\,du\le 0.14406,
\end{displaymath}
so that 
\begin{displaymath}
I_{12}<23.2227\frac{e^{-\pi \eta_2}}{|\eta|}.
\end{displaymath}

Hence,
\begin{equation}\label{firstR}
|R|\le28.9745\;
\frac{e^{-\pi \eta_2}}{|\eta|}+
14.2355\frac{e^{-\frac{\pi}{32}|\eta|^2}}{|\eta|}<29
\frac{e^{-\pi \eta_2}}{|\eta|}+ 15 \frac{e^{-\frac{\pi}{32}|\eta|^{2}}}{|\eta|}.
\end{equation}

Instead of the first integral on the right side of \eqref{split} we may also put (by \eqref{boundcos})
\begin{multline*}
I'_1:=\int_{|u|<|\eta|/4}\frac{e^{\pi/4} e^{-2\pi u^2}}{2^{-1/2} 
}\;\bigl(e^{\frac{4\pi}{3}\frac{|x-\eta|^3}{|\eta|}}-1\bigr)\,du\le \\ \le
2\sqrt{2}e^{\pi/4}\int_0^{|\alpha|}
e^{-2\pi u^2}\bigl(e^{\frac{4\pi}{3}\frac{2^{3/2}|\alpha|^3}{|\eta|}}-1\bigr)\,du+\\ +2\sqrt{2}e^{\pi/4}\int_0^{|\eta|/4}
e^{-2\pi u^2}\bigl(e^{\frac{4\pi}{3}\frac{2^{3/2}u^3}{|\eta|}}-1\bigr)\,du:=I'_{11}+I'_{1,2}.
\end{multline*}
We bound these two as before
\begin{displaymath}
I'_{11}\le \frac{2\pi\sqrt{2}e^{\pi/4}}{3|\eta|}\int_0^{2^{-3/2}}e^{-2\pi u^2}\,du
<1.8143\frac{1}{|\eta|}.
\end{displaymath}
\begin{displaymath}
I'_{12}\le 2\sqrt{2} e^{\pi/4} \frac{4\pi}{3}\frac{2^{3/2}}{|\eta|}\int_0^{|\eta|/4}
e^{-2\pi u^2}e^{\frac{4\pi}{3}\frac{2^{3/2}u^3}{|\eta|}} u^3\,du<2.0367\frac{1}{|\eta|}.
\end{displaymath}

So in addition to \eqref{firstR} we also have
\begin{equation}
|R|\le
3.8510\frac{1}{|\eta|}+
14.2355\frac{e^{-\frac{\pi}{32}|\eta|^2}}{|\eta|}<\frac{4}{|\eta|}+15\frac{e^{-\frac{\pi}{32}|\eta|^2}}{|\eta|}.
\end{equation}
So, combining these two we get our result. 
\end{proof}

\section{Main Terms of the approximation to \texorpdfstring{$\Rzeta(s)$}{R(s)}.}

\begin{lemma}\label{mainlemma}
In the closed set  $H$ defined by the inequalities
\begin{displaymath}
\eta_1\ge 2,\qquad \eta_2\ge2,
\end{displaymath}
we have 
\begin{equation}\label{defU1}
\Rzeta(s)=-\chi(s)\eta^{s-1}e^{-\pi i \eta^2}\frac{\sqrt{2}e^{3\pi i/8}\sin\pi\eta}
{2\cos2\pi\eta}(1+U),
\end{equation}
where 
\begin{multline}\label{initialbound}
|U|\le 1.10942\;\frac{\sqrt{2}}{\pi\eta_2}
(1+\pi\eta_2)e^{\pi\eta_2}\,e^{-2\pi\frac{\eta_1^2\eta_2^2}{\eta_1^2+\eta_2^2}} +\\+
\frac{\sqrt{2}}{\pi\eta_2}
(1+\pi\eta_2)e^{\pi\eta_2}\Bigl(
\frac{\min(4,29e^{-\pi \eta_2})}{|\eta|}+ 15\frac{e^{-\frac{\pi}{32}|\eta|^{2}}}
{|\eta|}\Bigr)+\frac{e^{-4\pi(\eta_2-\frac12)\eta_2}}{\sqrt{2}\sinh\pi\eta_2}.
\end{multline}
\end{lemma}

\begin{proof}
From the definition of $\eta=\eta_1+i\eta_2$ we get 
\begin{equation}\label{eta-s}
\eta_1^2-\eta_2^2=\frac{t}{2\pi},\qquad 2\eta_1\eta_2=\frac{1-\sigma}{2\pi}.
\end{equation}
In $H$ we have  $\frac{1-\sigma}{2\pi}=2\eta_1\eta_2>8$. Hence $\sigma<-10$, Theorem \ref{boundR} applies, and the zeta function $\zeta(1-s)$ is given by the Dirichlet series. Therefore, we can write \eqref{main} as
\begin{multline}\label{main21}
\Rzeta(s)=\chi(s)\eta^{s-1}e^{-\pi i \eta^2}\Bigl\{e^{\pi i \eta^2}\sum_{n=k+1}^{\infty}\Bigl(\frac{n}{\eta}\Bigr)^{s-1}-\\-
\frac{\sqrt{2}e^{3\pi i/8}\sin\pi\eta-(-1)^k e^{2\pi i\eta-2\pi i(\eta-k)^2}}
{2\cos2\pi\eta}-\\-\sum_{j=m+1}^k(-1)^{j-1}e^{-2\pi i(j-\eta)^2}
w(j-\eta)-R\Bigr\}.
\end{multline}

So, \eqref{defU1} is true if we define
\begin{multline}\label{defdefU}
U:=\frac{2\cos2\pi\eta}{\sqrt{2}e^{3\pi i/8}\sin\pi\eta}
\Bigl\{-e^{\pi i \eta^2}\sum_{n=k+1}^{\infty}\Bigl(\frac{n}{\eta}\Bigr)^{s-1}+\\
+\sum_{j=m+1}^k(-1)^{j-1}e^{-2\pi i(j-\eta)^2}
w(j-\eta)+R\Bigr\}-\frac{(-1)^k e^{2\pi i\eta-2\pi i(\eta-k)^2}}{\sqrt{2}e^{3\pi i/8}\sin\pi\eta}.
\end{multline}

Since $\frac{x\cosh2x}{\sinh x}<(1+x)e^x$ for $x>0$ we have
\begin{equation}\label{factor1}
\Bigl|-\frac{2\cos2\pi\eta}{\sqrt{2}e^{3\pi i/8}\sin\pi\eta}\Bigr|\le 
\frac{2\cosh2\pi\eta_2}{\sqrt{2}\sinh\pi\eta_2}<\frac{\sqrt{2}}{\pi\eta_2}
(1+\pi\eta_2)e^{\pi\eta_2}.
\end{equation}

To bound the first term in \eqref{defdefU} observe that since $\sigma<0$
\begin{multline}\label{boundfirstsum}
\Bigl|e^{\pi i \eta^2}\sum_{n=k+1}^{\infty}\Bigl(\frac{n}{\eta}\Bigr)^{s-1}\Bigr|<
e^{-2\pi\eta_1\eta_2}|\eta|^{1-\sigma}e^{t\arg\eta}\Bigl\{\frac{(k+1)^\sigma}{
-\sigma}+(k+1)^{\sigma-1}\Bigr\}<\\ <
e^{-2\pi\eta_1\eta_2+t\arg\eta}\Bigl(\frac{k+1}{|\eta|}\Bigr)^{\sigma-1}
\Bigl(\frac{k+1}{-\sigma}+1\Bigr).
\end{multline}

We shall use equation \eqref{main21} with $k=m$ so that the sum in $j$
in \eqref{defdefU} does not appear. (But in Theorem \ref{smallregion} we shall need equation \eqref{defdefU} for other values of $k$). 

Since $\eta_1>2$ and $\eta_2>2$ we have $0<\arg\eta<\frac{\pi}{2}$. 
When $t\le 0$ we have  $-2\pi\eta_1\eta_2+t\arg\eta<0$ and when $t>0$ we have
\begin{displaymath}
-2\pi\eta_1\eta_2+t\arg\eta<-2\pi\eta_1\eta_2+t\frac{\eta_2}{\eta_1}
=-2\pi\frac{\eta_2^3}{\eta_1}<0,
\end{displaymath}
in both cases  
\begin{equation}\label{I:some}
e^{-2\pi\eta_1\eta_2+t\arg\eta}<1.
\end{equation}

We have $1+2x>e^x$ for $0<x\le x_0\approx1.25643$. In particular, this inequality is true for 
$x=\frac{\eta_1\eta_2}{\eta_1^2+\eta_2^2}\le\frac12$. Since $\sigma<0$ and $\eta_1+\eta_2<m+1$ we get
\begin{displaymath}
\Bigl(\frac{m+1}{|\eta|}\Bigr)^{\sigma-1}<\Bigl(\frac{\eta_1^2+\eta_2^2+2\eta_1\eta_2}
{\eta_1^2+\eta_2^2}\Bigr)^{\frac{\sigma-1}{2}}<e^{\frac{\eta_1\eta_2}{\eta_1^2+\eta_2^2}
\frac{\sigma-1}{2}}=e^{-2\pi\frac{\eta_1^2\eta_2^2}{\eta_1^2+\eta_2^2}}
\end{displaymath}
so that
\begin{equation}
\Bigl|e^{\pi i \eta^2}\sum_{n=m+1}^{\infty}\Bigl(\frac{n}{\eta}\Bigr)^{s-1}\Bigr|
\le e^{-2\pi\frac{\eta_1^2\eta_2^2}{\eta_1^2+\eta_2^2}}\Bigl(1+\frac{1+\eta_1+\eta_2}{-\sigma}\Bigr).
\end{equation}

Since  $\eta_1+\eta_2>4$ we have $1+\eta_1+\eta_2<5(\eta_1+\eta_2)/4$ and since
 we assume $\sigma<-10$ we will have $-\sigma>\frac{10}{11}(1-\sigma)$ so  we get 
\begin{gather*}
\Bigl(1+\frac{1+\eta_1+\eta_2}{-\sigma}\Bigr)
\le 1+\frac{11}{8}\frac{\eta_1+\eta_2}{1-\sigma}=
1+\frac{11}{8}\frac{\eta_1+\eta_2}{4\pi\eta_1\eta_2}=\\=
1+\frac{11}{32\pi}(\eta_1^{-1}+\eta_2^{-1})\le 1+\frac{11}{32\pi}<1.10942
\end{gather*}
so that
\begin{equation}\label{ineqsum1}
\Bigl|e^{\pi i \eta^2}\sum_{n=m+1}^{\infty}\Bigl(\frac{n}{\eta}\Bigr)^{s-1}\Bigr|
\le 1.10942\, e^{-2\pi\frac{\eta_1^2\eta_2^2}{\eta_1^2+\eta_2^2}}.
\end{equation}

By the definition of $m$ we have $m\le \eta_1+\eta_2<m+1$ so that
\begin{equation}\label{lastterm1}
|(-1)^{m-1}e^{2\pi i\eta-2\pi i(\eta-m)^2}|=e^{-2\pi\eta_2-4\pi(m-\eta_1)\eta_2}<e^{-4\pi(\eta_2-\frac12)\eta_2}.
\end{equation}

Joining the definition of $U$ in \eqref{defdefU} with  $k=m$  and \eqref{factor1}, \eqref{ineqsum1}, \eqref{lastterm1} and the bound of $R$ \eqref{boundReq1} we get the conclusion \eqref{initialbound} of our lemma. 
\end{proof}

\begin{theorem}\label{mainTheorem}
In the closed set  $G$ defined by the inequalities
\begin{displaymath}
\sigma<0,\quad |s-1|\ge 5408\pi,\quad t\ge-\frac{99}{20}(1-\sigma),\quad 
(1-\sigma)^2\ge225\pi t,
\end{displaymath}
we have 
\begin{equation}\label{defU}
\Rzeta(s)=-\chi(s)\eta^{s-1}e^{-\pi i \eta^2}\frac{\sqrt{2}e^{3\pi i/8}\sin\pi\eta}
{2\cos2\pi\eta}(1+U),
\end{equation}
where $|U|<0.9$.

It follows that the only zeros on this closed set are the trivial zeros at even negative integers. It also follows that they are simple zeros.
\end{theorem}

\begin{figure} \centering
  \includegraphics[width=9cm]{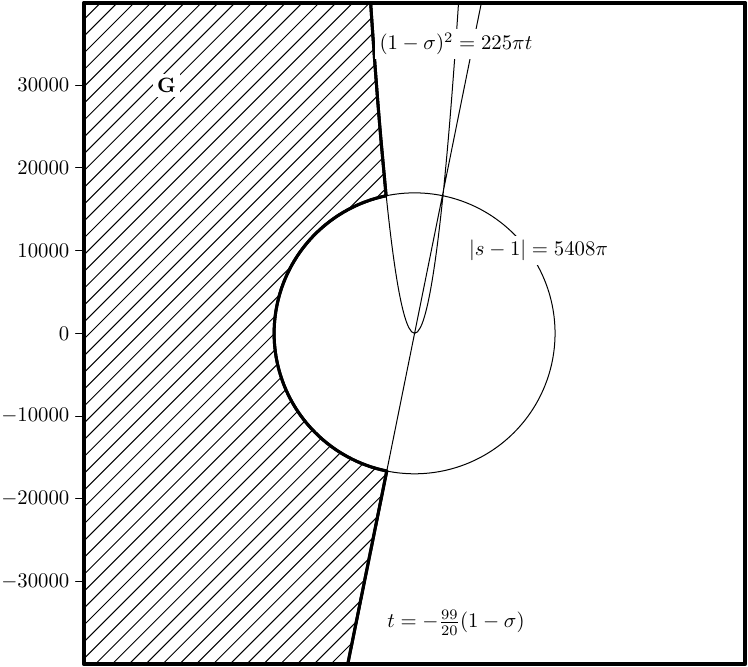}\\
  \caption{Region G of Theorem \ref{mainTheorem}}\label{figure1}
\end{figure}

\begin{proof}
Figure \ref{figure1} represents the region $G$. Since $\sigma<0$ the point $\eta$ is on the first quadrant, therefore $\eta_1$, $\eta_2>0$. The inequality $|s-1|\ge5408\pi$ is equivalent to $|\eta|\ge52$.  The inequality $t\ge-\frac{99}{20}(1-\sigma)$ restrict 
$\frac{s-1}{2\pi i}$ to be in a half-plane, so it is equivalent to $-\alpha<\arg\eta<\frac{\pi}{2}-\alpha$ where $\tan\alpha=\tan(\frac12\arctan\frac{20}{99})=\frac{1}{10}$. It follows that the inequality (~for $\eta_1>0$ and $\eta_2>0$~) is equivalent  to $\eta_2\le a\,\eta_1$ with $a=10$. 

Finally, the last  inequality $(1-\sigma)^2\ge225\pi t$ implies that $\eta_2\ge 5$. In fact by equations \eqref{eta-s} we  have the following
\begin{displaymath}
4\eta_1^2\eta_2^2\ge \frac{225}{2}(\eta_1^2-\eta_2^2),\quad\text{so that}\quad
\Bigl(4\eta_1^2+\frac{225}{2}\Bigr)\eta_2^2\ge\frac{225}{2}\eta_1^2.
\end{displaymath}
Now, if we assume that $\eta_2<5$ we get 
\begin{displaymath}
4\eta_1^2+\frac{225}{2}>\frac{9}{2}\eta_1^2
\end{displaymath}
equivalent to $\eta_1^2<225$ and we get the contradiction $52^2\le |\eta|^2\le \eta_1^2+\eta_2^2<250$.

Hence, for $s\in G$ we have
\begin{equation}\label{region}
|\eta|\ge 52,\quad 5\le \eta_2\le a\,\eta_1, \qquad a=10.
\end{equation}
In what follows, we assume that $\eta_1$ and $\eta_2$ satisfies these inequalities.
\medskip

Observe that Lemma \ref{mainlemma} applies, because for $\eta_2>20$ we have $\eta_1\ge \frac{\eta_2}{a}\ge2$, and for $\eta_2<20$, we have $52^2\le \eta_1^2+\eta_2^2$, so that $\eta_1\ge 
\sqrt{52^2-20^2}=48$. 
\medskip 

When $\eta_2\le \eta_1$ we have $2\pi\frac{\eta_1^2\eta_2^2}{\eta_1^2+\eta_2^2}\ge\pi\eta_2^2$
and (being $\eta_2>5$)
\begin{displaymath}
1.10942\,\frac{\sqrt{2}}{\pi\eta_2}
(1+\pi\eta_2)e^{\pi\eta_2}e^{-2\pi\frac{\eta_1^2\eta_2^2}{\eta_1^2+\eta_2^2}}\le 2\,\frac{\sqrt{2}}{\pi\eta_2}
(1+\pi\eta_2)e^{\pi\eta_2-\pi\eta_2^2}<2\times 10^{-27}.
\end{displaymath}
when $\eta_2>\eta_1$ we have $\eta_1=|\eta|\cos\varphi$, $\eta_2=|\eta|\sin\varphi$ with $\frac{\pi}{4}\le\varphi\le \arctan 10$. Then 
\begin{multline*}
\pi\eta_2-2\pi\frac{\eta_1^2\eta_2^2}{\eta_1^2+\eta_2^2}
=\pi |\eta|(\sin\varphi-2|\eta|\sin^2\varphi\cos^2\varphi)\le \\ \le
\pi |\eta|(\sin\varphi-104\sin^2\varphi\cos^2\varphi)<-0.076876 |\eta|,
\end{multline*}
also, it is clear that $\pi \eta_2=\pi|\eta|\sin\varphi>52\pi/\sqrt{2}:=b$, so that
\[
1.10942\frac{\sqrt{2}}{\pi\eta_2}
(1+\pi\eta_2)e^{\pi\eta_2}e^{-2\pi\frac{\eta_1^2\eta_2^2}{\eta_1^2+\eta_2^2}} \le1.10942\frac{\sqrt{2}}{b}
(1+b)e^{-0.076876 |\eta|}<0.0290564.
\]
For $\eta_2\ge5$ we have $4>29e^{-\pi\eta_2}$ and
\begin{displaymath}
\frac{\sqrt{2}}{\pi\eta_2}
(1+\pi\eta_2)e^{\pi\eta_2}\frac{29e^{-\pi\eta_2}}{|\eta|}\le \frac{\sqrt{2}}{\pi\eta_2}
(1+\pi\eta_2)\frac{29}{\max(52, \eta_2)}<0.83891.
\end{displaymath}
(This is where we need $52$ so large, and  I would very much like a smaller radius.)
For $|\eta|>52$ and $\eta_2>5$ 
\begin{displaymath}
\frac{\sqrt{2}}{\pi\eta_2}
(1+\pi\eta_2)e^{\pi\eta_2}\frac{15 e^{-\frac{\pi}{32}|\eta|^2}}{|\eta|}\le 
\frac{\sqrt{2}}{5 \pi}
(1+5\pi)e^{\pi|\eta|}\frac{15 e^{-\frac{\pi}{32}|\eta|^2}}{|\eta|}<3\times10^{-45}
\end{displaymath}
and for $\eta_2>5$ 
\begin{displaymath}
\frac{1}{\sqrt{2}\sinh\pi\eta_2}e^{-4\pi(\eta_2-\frac12)\eta_2}<4\times 10^{-130}.
\end{displaymath}

Collecting these inequalities, we get 
\begin{displaymath}
|U|\le 0.0290564+0.83891+3\times10^{-45}+4\times 10^{-130}<0.868.
\end{displaymath}

The assertion about the zeros follows easily from the representation \eqref{defU}.
\end{proof}

\begin{theorem} \label{apRzeta}
Let  $G_\alpha$ be the closed set defined by  $|t|\cdot\tan2\alpha\le (1-\sigma)$, where $0<\alpha<\frac{\pi}{4}$.  For $s\to\infty$ we have
\begin{equation}\label{apRzetaU}
\Rzeta(s)=\frac{e^{-\frac{\pi i}{8}}}{\sqrt{2}}\chi(s)\eta^{s-1}e^{\pi i(\eta- \eta^2)}(1+\Orden(|\eta|^{-1})),\qquad  s\in G_\alpha.
\end{equation}
\end{theorem}

\begin{proof}

The assertion $s\in G_\alpha$ is equivalent to $2\alpha\le\arg(\frac{s-1}{2\pi i})\le \pi-2\alpha$. Therefore, it is equivalent  to $\alpha<\arg(\eta)<\frac{\pi}{2}-\alpha$. 
Therefore, when  $s\in G_\alpha$ we have $\eta=\eta_1+i\eta_2$ with $\eta_j\ge0$. In fact  we have 
$|\eta|\sin\alpha<\eta_j<|\eta|\cos\alpha$ for $j=1$ and  $j=2$.
Since $\alpha>0$, the two variables $\eta_j\to\infty$ when $s$ goes to infinity. So, to prove our Theorem, we may assume that $\eta_j\ge2$, and Lemma \ref{mainlemma} applies.

Therefore,
\[
\pi\eta_2-2\pi\frac{\eta_1^2\eta_2^2}{\eta_1^2+\eta_2^2}<\pi|\eta|\cos\alpha-2\pi
\frac{|\eta|^4\sin^4\alpha}{|\eta|^2}<-(\pi\cos\alpha)|\eta|,\quad
\text{for}\quad |\eta|>\frac{\cos\alpha}{\sin^4\alpha}.
\]
For $|\eta|>\sqrt{2}$, we have $\eta_2>
\sqrt{2}\sin\alpha$. Hence, for $|s|$ large enough
\begin{displaymath}
1.10942\frac{\sqrt{2}}{\pi\eta_2}
(1+\pi\eta_2)e^{\pi\eta_2}e^{-2\pi\frac{\eta_1^2\eta_2^2}{\eta_1^2+\eta_2^2}}=\Orden(e^{-(\pi\cos\alpha)|\eta|}),\qquad |\eta|\to\infty
\end{displaymath}
It follows easily from \eqref{initialbound} that $|U|=\Orden(|\eta|^{-1})$.

Finally, we observe that for $s\in G_\alpha$ when $s\to\infty$ we have
\begin{displaymath}
\frac{e^{3\pi i/8}\sin\pi\eta}{\cos2\pi\eta}= \frac{-e^{3\pi i/8}e^{-\pi i\eta}}{i e^{-2\pi i\eta}}(1+\Orden(|\eta|^{-1}))=e^{7\pi i/8} e^{\pi i\eta}(1+\Orden(|\eta|^{-1})).
\end{displaymath}
\end{proof}

\begin{corollary}
The only zeros of $\Rzeta(s)$ in the region defined by the inequalities $1-\sigma>|t|$, $|s-1|\ge3528\pi$  are the trivial zeros at the points $s=-2n$ for $n$ a positive integer.
\end{corollary}

\begin{proof}
As in  the proof of  Theorem \ref{apRzeta} we take $\alpha=\pi/8$ so that $\tan2\alpha=1$.
For $|\eta|\ge42$ we have $\eta_j>|\eta|\sin\alpha>2$ and Lemma \ref{mainlemma} applies.
We have 
\begin{multline*}
1.10942\,\frac{\sqrt{2}}{\pi\eta_2}
(1+\pi\eta_2)e^{\pi\eta_2}e^{-2\pi\frac{\eta_1^2\eta_2^2}{\eta_1^2+\eta_2^2}}\le\\ \le 1.10942\,\sqrt{2}\Bigl(1+\frac{1}{\pi|\eta|\sin\alpha}\Bigr)
e^{\pi|\eta|\cos\alpha-2\pi|\eta|^2\sin^4\alpha}<
9\times10^{-50}.
\end{multline*}
For $|\eta|\ge42$ we have $29 e^{-\pi\eta_2}<29e^{-\pi |\eta|\sin\alpha}<4$ and
\begin{multline*}
\frac{\sqrt{2}}{\pi\eta_2}
(1+\pi\eta_2)e^{\pi\eta_2}\Bigl(
\frac{\min(4,29e^{-\pi\eta_2})}{|\eta|}+ 15\frac{e^{-\frac{\pi}{32}|\eta|^{2}}}{|\eta|}\Bigr)\le\\ \le \sqrt{2}\Bigl(1+\frac{1}{\pi|\eta|\sin\alpha}\Bigr)
\Bigl(\frac{29}{|\eta|}+15\frac{e^{\pi |\eta|\cos\alpha-
\pi|\eta|^2/32}}{|\eta|}\Bigr)\le 0.99582
\end{multline*}
and finally for $|\eta|>42$, we have $\eta_2>42\sin\alpha$ and 
\begin{displaymath}
\frac{1}{\sqrt{2}\sinh\pi\eta_2}e^{-4\pi(\eta_2-\frac12)\eta_2}\le 
2\times 10^{-1388}.
\end{displaymath}
By \eqref{initialbound} we get $|U|<1$ on our closed set. Then \eqref{apRzetaU} implies  that the zeros of $\Rzeta(s)$ in the region 
limited by $1-\sigma>|t|$ and $|\eta|>42$ are those of $\chi(s)$. These are the trivial zeros at  even negative integers. 
\end{proof}

\begin{remark}
Computations of the zeros of $\Rzeta(s)$ indicate that 
the zeros in the set  of points in the  circle $|\eta|\le42$ and on the angle $1-\sigma>|t|$ are just the trivial ones.  
\end{remark}

\section{Left limit of the zeros for \texorpdfstring{$t>0$}{t>0}.}

From Theorem \ref{mainTheorem} we know that in the region defined by $t>0$,  $|s-1|>2\pi b^2$, $1-\sigma\ge a t^{1/2}$ (where $b=52$ and $a=15\sqrt{\pi}$) there is no zero of $\Rzeta(s)$. Now we are going to extend this region, changing the condition $1-\sigma\ge a t^{1/2}$ for the less restrictive $1-\sigma\ge c t^{2/5}\log t$. Siegel says in \cite{Siegel} that given $\varepsilon>0$ he may prove this for the region  $1-\sigma\ge c t^{\varepsilon}$, but I have not been able (following his method) to get a better exponent than $2/5$.  Siegel only proves the case of  exponent $3/7$. 

\begin{theorem}\label{T:once}
There are constants $A>0$ and $t_0>0$ such that 
\begin{equation}\label{defU2}
\Rzeta(s)=-\chi(s)\eta^{s-1}e^{-\pi i \eta^2}\frac{\sqrt{2}e^{3\pi i/8}\sin\pi\eta}
{2\cos2\pi\eta}(1+U),\qquad |U|<1,
\end{equation}
for $t>t_0$ and $1-\sigma\ge A\; t^{2/5}\log t$.
\end{theorem}

\begin{proof} By Theorem \ref{mainTheorem} we know that \eqref{defU2} is true for 
$t>t_0$ and $1-\sigma\ge a t^{1/2}$. Therefore, we shall consider the region $\Omega$ defined by $t>t_0$ and $1<\varphi(t)\le 1-\sigma\le a t^{1/2}$.  Our objective is to get $\varphi(t)$ as small as possible.  The proof is similar to the one of Lemma \ref{mainlemma}. For $s\in\Omega$ the hypothesis of Theorem \ref{boundR} holds, and we have 
\eqref{main} with $R$ bounded as in \eqref{boundReq1}. When $s\in\Omega$ we have $\sigma_0<-10$ so that \eqref{main} can be written as \eqref{main21}. 
We take  $k=m+r$, and  $r$ will be a function of $s=\sigma+it$, $r=\lfloor \psi\rfloor$. (Siegel uses $\psi(\sigma)$ but we consider a more general $\psi$ that may be a function of $\sigma$ and $t$).
We will use the letter $c$ to denote a constant, not always the same, in particular, note that $r=\lfloor \psi\rfloor$, but we will choose  $\psi\to+\infty$ for $t\to+\infty$ so that $ c\psi \le r\le \psi$.

Our elections will be $\varphi(t)=A\;t^{2/5}\log t$ and $\psi=\psi(t)=Bt^{1/10}\log^{1/5} t$ for adequate constants $A$ and $B$. We retain in our reasoning a general $\varphi(t)$ and $\psi$ to see that the reasoning does not give a better (~smaller~) $\varphi(t)$.

We have $1-\sigma=\Orden(t^{1/2})$ for $t\to\infty$, therefore for 
$s\in\Omega$ and $t>t_0$ we have the convergent expansion
\begin{equation}\label{etat}
\eta=\Bigl(\frac{t}{2\pi}+\frac{1-\sigma}{2\pi}i\Bigr)^{1/2}=
\sqrt{\frac{t}{2\pi}}\Bigl(1+i\frac{1-\sigma}{2t}+\frac{(1-\sigma)^2}{8t^2}-i\frac{(1-\sigma)^3}{16t^3}-\frac{(1-\sigma)^4}{128t^4}+\cdots\Bigr)
\end{equation}
so that
\begin{equation}\label{eta12t}
\eta=\Bigl(\frac{t}{2\pi}\Bigr)^{1/2}+\Orden(1),\quad
\eta_1=\Bigl(\frac{t}{2\pi}\Bigr)^{1/2}+\Orden(t^{-1/2}),\quad \eta_2=\frac{1-\sigma}{2\sqrt{2\pi t}}+\Orden(t^{-1}).
\end{equation}

Now we must bound $U$ given in \eqref{defU1} and \eqref{defdefU}. Observe that \eqref{factor1} and 
\eqref{boundfirstsum} are still  true in this context. 

With respect to \eqref{boundfirstsum} observe that $m\le \eta_1+\eta_2<m+1$ and $|\eta|\le \eta_1+\eta_2<m+1$ so that $\frac{|\eta|+r}{|\eta|}<\frac{m+r+1}{|\eta|}$ and therefore for $t>t_0$ 
(~notice that $1+x\ge e^{2x/3}$ for $0\le x\le1$, and take $t_0$ such that for $t>t_0$
$r<|\eta|<\frac43\eta_1$. By \eqref{eta12t} we have $|\eta|\sim\eta_1$ for $t\to\infty$.)
\[
\Bigl(\frac{m+r+1}{|\eta|}\Bigr)^{\sigma-1}<\Bigl(\frac{|\eta|+r}{|\eta|}\Bigr)^{\sigma
-1}<e^{(\sigma-1)\frac{2r}{3|\eta|}}=e^{-(1-\sigma)\frac{2r}{3|\eta|}}<2
e^{-(1-\sigma)\frac{r}{2\eta_1}}=2e^{-2\pi r\eta_2}.
\]
The assumption $\frac{r}{|\eta|}<1$,  will follow from our final choice of $\psi$ (that will be $\Orden(t^{1/10}\log^{1/5}t)$) and \eqref{etat}.
It follows that for $t\ge t_0$ and $s\in\Omega$
\begin{displaymath}
\Bigl(\frac{m+r+1}{|\eta|}\Bigr)^{\sigma-1}\le 2 e^{-c(1-\sigma) \psi t^{-1/2}}.  
\end{displaymath}

For $s\in\Omega$ we have $\varphi(t) t^{-1/2}\le (1-\sigma)t^{-1/2}\le a$, and $\eta_2$ is given by \eqref{eta12t}, so that for $t>t_0$ we have
$c_1\frac{\varphi(t)}{t^{1/2}}\le \eta_2\le c_2$. It follows that 
\begin{displaymath}
\frac{\sqrt{2}}{\pi\eta_2}
(1+\pi\eta_2)e^{\pi\eta_2}\le c\frac{t^{1/2}}{\varphi(t)}
\end{displaymath}
so that (notice that \eqref{I:some} holds in this case)
\begin{multline*}
\frac{\sqrt{2}}{\pi\eta_2}
(1+\pi\eta_2)e^{\pi\eta_2}\Bigl(\frac{m+r+1}{|\eta|}\Bigr)^{\sigma-1}\Bigl(1+\frac{m+r+1}{-\sigma}\Bigr)\le\\ \le c\frac{t^{1/2}}{\varphi(t)}e^{-c(1-\sigma)\psi t^{-1/2}}\Bigl(1+c \frac{t^{1/2}+\psi}{\varphi(t)}\Bigr)\le c \frac{te^{-c(1-\sigma)\psi t^{-1/2}}}{\varphi(t)^{2} },
\end{multline*}
because we shall take $\psi\le t^{1/2}$ and $\varphi(t)=A\;t^{2/5}\log t$. 

Now we consider other terms  in \eqref{defdefU}. First, for $1\le\ell\le r$, since $\eta_2 =\Orden(1)$ we have
\begin{displaymath}
|m+\ell-\eta|^2=(m+\ell-\eta_1)^2+\eta_2^2\le (r+\eta_2)^2+\eta_2^2=\Orden(\psi^2).
\end{displaymath}
We choose $\psi(\sigma)=\orden(t^{1/2})$ so that  $|m+\ell-\eta|<\frac12 |\eta|$ 
for $t\ge t_0$. Then we may apply \eqref{restrineq} and
\begin{equation}
|w(m+\ell-\eta)|\le e^{c\frac{\psi^3}{|\eta|}}-1\le c\frac{\psi^3}{t^{1/2}},
\end{equation}
(again, by our final election of $\psi$).

Also,
\begin{displaymath}
|e^{-2\pi i(m+\ell-\eta)^2}|= e^{-4\pi(m+\ell-\eta_1)\eta_2}\le  1,
\end{displaymath}
since $m+\ell-\eta_1\ge m+1-\eta_1>\eta_2$. 
Then
\begin{displaymath}
\frac{\sqrt{2}}{\pi\eta_2}
(1+\pi\eta_2)e^{\pi\eta_2}\Bigl|\sum_{j=m+1}^{m+r}(-1)^{j-1}e^{-2\pi i(j-\eta)^2}
w(j-\eta)\Bigr|\le cr \frac{t^{1/2}}{\varphi(t)}\frac{\psi^3}{t^{1/2}}\le 
c\frac{\psi^4}{\varphi(t)}.
\end{displaymath}

For $t>t_0$ we will also have
\begin{displaymath}
|e^{2\pi i\eta-2\pi i(m+r-\eta)^2}|=e^{-2\pi \eta_2-4\pi(m+r-\eta_1)\eta_2}<e^{-2\pi \eta_2-4\pi (r-1+\eta_2)\eta_2}<e^{-4\pi (r-1) \eta_2}.
\end{displaymath}
Hence, for $t>t_0$ we have
\[
\Bigl|\frac{(-1)^k e^{2\pi i\eta-2\pi i(\eta-k)^2}}{\sqrt{2}e^{3\pi i/8}\sin\pi\eta}\Bigr|\le 
\frac{e^{-4\pi (r-1)\eta_2}}{\sqrt{2}\sinh\pi\eta_2}\le c\frac{e^{-c(1-\sigma)\psi t^{-1/2}}}{\pi\eta_2}
\le c\frac{t^{1/2}}{\varphi(t)}e^{-c(1-\sigma)\psi t^{-1/2}}.
\]

The rest $R$ in \eqref{defdefU} by Theorem \ref{boundR} is bounded by
\begin{displaymath}
\frac{\sqrt{2}}{\pi\eta_2}
(1+\pi\eta_2)e^{\pi\eta_2}|R|\le c\frac{t^{\frac12}}{\varphi(t)}\frac{c}{|\eta|}\le 
\frac{c}{\varphi(t)}.
\end{displaymath}

In summary we have
\begin{displaymath}
|U|\le c t \varphi(t)^{-2}e^{-c(1-\sigma)\psi t^{-1/2}}+
c\psi^4\varphi(t)^{-1}+\frac{c}{\varphi(t)}+ct^{1/2}\varphi(t)^{-1} e^{-c(1-\sigma)\psi 
t^{-1/2}}.
\end{displaymath}
We will not expect this to be $<1$ for $t\to+\infty$. The first terms forces   $(1-\sigma)\psi(\sigma,t)\gg t^{1/2}\log t$ and the second $\psi^4\ll \varphi(t)$. Hence,
\begin{displaymath}
\frac{t^{1/2}\log t}{1-\sigma}\ll \psi\ll \varphi(t)^{1/4}
\end{displaymath}
This must be true for all values of $\sigma$. In particular, we must have
\begin{displaymath}
\frac{t^{1/2}\log t}{\varphi(t)}\ll \varphi(t)^{1/4}. 
\end{displaymath}
That is,
\begin{displaymath}
\varphi(t)\gg t^{2/5}\log^{4/5}t.
\end{displaymath}

This shows that in this way we may not get $\varphi(t)\sim t^\varepsilon $ for any $\varepsilon<2/5$.

If we choose $\varphi(t)= A t^{2/5}\log t$ and 
$\psi(t)=B t^{1/10}\log^{1/5}t$, with an adequate election of $A$ and $B$,  we get $U=\Orden(\log^{-1/5}t)$.
\end{proof}

\begin{corollary}\label{smallregion}
There exist constants $t_0$ and $A$ such that 
 $\Rzeta(s)$ have no zeros in the region limited by 
$t\ge t_0$, $A\, t^{2/5}\log t> 1-\sigma $.
\end{corollary}


\begin{thebibliography}{999}

\bibitem{A86}
\textsc{J. Arias de Reyna},  \href{https://doi.org/10.1090/S0025-5718-2010-02426-3 }{\emph{High precision computation of Riemann's zeta function by the Riemann-Siegel formula, I}},  Math. Comp. \textbf{80} (2011), no. 274, 995--1009.

\bibitem{A166}
\textsc{Arias de Reyna, J.},  \emph{Riemann's auxiliary function. Basic results}, preprint. \href{https://arxiv.org/abs/2406.02403}{arXiv:2406.02403} 


\bibitem{A172}
\textsc{Arias de Reyna, J.},  \emph{Statistics of zeros of the auxiliary function}, preprint. \href{https://arxiv.org/abs/2406.03041}{arXiv:2406.03041} 

\bibitem{G} \textsc{W.~Gabcke}, \emph{Neue Herleitung und
explizite Restabschätzung der Rie\-mann\-{}-Siegel\-{}-Formel},
Dissertation, Göttingen (1979).  \href{https://ediss.uni-goettingen.de/bitstream/handle/11858/00-1735-0000-0022-6013-8/DissGab1979.pdf?sequence=1}{Electronic version.}


\bibitem{O'S}\textsc{Cormac O'Sullivan}, \emph{A Generalization of the Riemann-Siegel formula}, Math. Z. \textbf{303} 20 (2023). \href{https://arxiv.org/abs/1811.01130}{arXiv:1811.01130}.

\bibitem{Siegel}
\textsc{C. L. Siegel}, \emph{Über Riemann Nachla\ss\ zur analytischen
Zahlentheorie}, Quellen und Studien zur Geschichte der Mathematik Astronomie und 
Physik \textbf{2} (1932) 45--80. (Reprinted in \cite{SW}, 1, 275--310.)
\href{https://arxiv.org/abs/1810.05198}{English translation}.

\bibitem{SW}
\textsc{C. L. Siegel}, \emph{Carl Ludwig Siegel's Gesammelte Abhandlungen}, 
(edited by K. Chandrasekharan and H. Maa\ss), Springer-Verlag, Berlin, 1966.


\end{thebibliography}
\end{document}